\documentclass[12pt]{amsart}

\usepackage{mathrsfs}
\usepackage{amsmath}
\usepackage{latexsym}
\usepackage{amssymb}
\usepackage{amsfonts}
\usepackage{cite}

\include{PDF}
\usepackage{graphics}
\def\l{\langle} \def\r{\rangle} 
\def\div{\,\,\big|\,\,} 
 \def\ZZ{\mathbb Z}

 \def\val{{\sf val}}

\def\Aut{{\sf Aut}} 

\def\Cos{{\sf Cos}}

\def\Cay{{\sf Cay}} 
\def\D{{\rm D}} 
\def\S{{\rm S}} 
 \def\M{{\rm M}}

\def\C{{\rm C}}\def\N{{\rm N}} \def\O{{\bf O}}

\def\Ga{{\it \Gamma}}

\def\a{\alpha}   \def\s{\sigma}

 \def\GL{{\rm GL}}

\def\GammaL{{\rm \Gamma L}}
\def\AGammaL{{\rm A\Gamma L}}\def\ASigmaL{{\rm A\Sigma L}}
 
\def\A{{\sf A}}

\def\PSL{{\rm PSL}}
\def\GL{{\rm GL}} 
\def\AGL{{\rm AGL}}\def\ASL{{\rm ASL}}

 \def\F{{\rm F}} \def\D{{\rm D}}

\newtheorem{rmk}{Remark}[section]
\newtheorem{thm}{Theorem}[section]

\newtheorem{construct}[thm]{Construction}

\newtheorem{lemma}{Lemma}[section]

\def\pf{\noindent{\it Proof.} }
\def\qed{\nopagebreak\hfill{\rule{4pt}{7pt}}
\medbreak}

\begin{document}

\title[Cayley graphs on simple groups]{Arc-transitive pentavalent Cayley graphs with soluble vertex stabilizer on finite nonabelian simple groups$^{*}$}
\thanks{2000 MR Subject Classification 20B25, 05E18, 05C25.}
\thanks{$^{*}$This work was partially supported by the National Natural Science Foundation of China (11231008,11461004),
and the Natural Science Foundation of Yunnan Province (2013FB001, 2015J006).}
\author[B. Ling and B.G. Lou] {Bo Ling$^{1}$ and Ben Gong Lou$^{2}$}
\address{1: School of Mathematics and Computer Sciences\\
Yunnan Minzu University\\
Kunming, Yunnan 650031, P.R. China}
\email{bolinggxu@163.com (B. Ling)}
\address{2: School of Mathematics and Statistics\\ Yunnan University,
Kunmin 650031, P. R. China}
\email{bengong188@163.com (B.G. Lou)}

\begin{abstract}
A Cayley graph $\Ga=\Cay(G,S)$ is said to be normal if $G$ is normal in $\Aut\Ga$.
The concept of normal Cayley graphs was first proposed by M.Y.Xu in [Discrete Math. 182, 309-319, 1998] and
it plays an important role in determining the full automorphism groups of Cayley graphs.
In this paper, we investigate the normality problem of the connected arc-transitive pentavalent Cayley graphs
with soluble vertex stabilizer on finite nonabelian simple groups. We prove that all such graphs $\Ga$ are either normal or
$G=\A_{39}$ or $\A_{79}$. Further, a connected arc-transitive pentavalent Cayley graph on $\A_{79}$ is constructed.
To our knowledge, this is the first known example of pentavalent 3-arc-transitive Cayley graph
on finite nonabelian simple group which is non-normal.
\vskip4pt

\noindent {\sc Keywords}. Simple group; Normal Cayley graph; Arc-transitive graph
\end{abstract}
\maketitle

\section{ Introduction}
All graphs are assumed to be finite, simple and undirected.

Let $\Ga$ be a graph. We use $V\Ga$, $E\Ga$
and $\Aut\Ga$ to denote the vertex set, edge set and
automorphism group of $\Ga$, respectively.
Denote $\val\Ga$ the valency of
$\Ga$. Let $X\le\Aut\Ga$ and
let $s$ be a positive integer. The graph $\Ga$ is said to be {\em $(X,s)$-arc-transitive},
if $X$ acts transitively on the set of $s$-arcs of $\Ga$, where an {\em $s$-arc} is an $(s+1)$-tuple $(v_0, v_1,\cdots,v_s)$
of $s+1$ vertices satisfying $(v_{i-1},v_i)\in E\Ga$ and $v_{i-1}\not=v_{i+1}$ for all $i$.
The graph $\Ga$ is called {\em $(X,s)$-transitive} if it is $(X,s)$-arc-transitive
but not $(X,s+1)$-arc-transitive. In particular, an $(\Aut\Ga,s)$-arc-transitive
or $(\Aut\Ga,s)$-transitive graph is just called {\em $s$-arc-transitive}
or {\em $s$-transitive graph}; and $0$-arc-transitive
graph is called {\em vertex transitive} graph, $1$-arc-transitive graph is called
{\em arc-transitive} graph or {\em symmetric} graph.

Let $G$ be a finite group with identity $1$, and let $S$ be a subset of $G$ such that
$1\not\in S$ and $S=S^{-1}:=\{x^{-1}\mid x\in S\}$. The  Cayley graph of $G$ with respect to $S$, denoted by $\Cay(G,S)$,
is defined on $G$ such that $g,\,h\in G$ are adjacent if and only if $hg^{-1}\in S$.
Then $\Cay(G,S)$ is a regular graph of valency $|S|$. It is well-known that $\Ga$ is connected if and only if $\l S\r=G$, that is,
$S$ is a generating set of the group $G$. For a Cayley graph $\Cay(G,S)$, the underlying group $G$ can be viewed as a regular subgroup of $\Aut\Cay(G,S)$ which acts on $G$ by right multiplication.
Conversely, a graph $\Ga$ is isomorphic to a Cayley graph if and only if $\Aut\Ga$ has a regular subgroup, refer to \cite[Proposition 16.3]{Biggs}.

Let $\Ga=\Cay(G,S)$ be a Cayley graph. Set
$$\Aut(G,S)=\{\s\in \Aut(G)\mid S^\s=S\}.$$
Then $\Aut(G,S)$, acting on $G$ naturally, is a subgroup of $\Aut\Ga$. If $\Ga$ is connected, then $\Aut(G,S)$ acts faithfully on $S$ and lies in
the stabilizer of the vertex corresponding to the identity of $G$.
Moreover, the normalizer
$\N_{\Aut\Ga}(G)$ equals to the semi-directed product $G{:}\Aut(G,S)$, see \cite{Godsil}.

A Cayley graph $\Ga=\Cay(G,S)$ is said to be {\em normal} if $G$ is normal in $\Aut\Ga$, that is, $\Aut\Ga= G{:}\Aut(G,S)$,
refer to \cite{Xu98}; otherwise, $\Ga$ is called {\em non-normal}. Thus, for a connected normal Cayley graph $\Ga=\Cay(G,S)$,
the group $\Aut(G,S)$ is just the stabilizer in $\Aut\Ga$ of the vertex corresponding to the identity of $G$.

In this paper we consider connected arc-transitive pentavalent Cayley graphs.

The concept of normal Cayley graphs was first proposed by M.Y.Xu in \cite{Xu98} and
it plays an important role in determining the full automorphism groups of Cayley graphs.
The Cayley graphs on finite nonabelian simple groups are received most attention in the
literature. For example, X.G.Fang, C.E.Praeger and J.Wang \cite{FPW02} gave a general description of the
possibilities for the automorphism groups of connected Cayley graphs on a finite non-abelian
simple group. Then some further work focuses on the small valencies because the precise structure
of the vertex stabilizer of arc-transitive cubic, tetravalent and pentavalent graphs was
determined by a series of papers.
Let $\Ga=\Cay(G,S)$ be a connected arc-transitive Cayley graph on a
finite nonabelian simple group $G$. For the cubic case, C.H.Li \cite{Li96} proved that only $7$ groups
are exceptions for $\Ga$ being not normal; on the basis of C.H.Li's result,
S.J.Xu et al. \cite{Xufang,XufangWang} proved that all such $\Ga$ are normal except two
$5$-transitive Cayley graphs of the alternating group $\A_{47}$, and so a complete classification
of cubic $s$-transitive non-normal Cayley graphs of finite simple groups was given.
For the tetravalent case, X.G.Fang, C.H.Li and M.Y.Xu in \cite{FangLiXu-01} proved that
most of such graphs are normal except a list of possible $G$.
Recently, J.J.Li, J.C.Ma and the first author of the present paper in \cite{LiLing01}
proved that all $s$-regular $\Ga$ (that is, $\Aut\Ga$ acts regularly on its $s$-arc set) are
normal except the case $G=\A_{35}$, and a 3-arc-transitive non-normal Cayley graph on $\A_{35}$ was constructed.
Further, X.G.Fang, J.Wang and S.M.Zhou in \cite{FWZ16} proved that all $2$-transitive $\Ga$ are normal except two graphs on $\M_{11}$.
For the pentavalent case, J.X.Zhou and Y.Q.Feng in \cite{Zhou-1} proved that all 1-transitive $\Ga$ are normal,
and in \cite{LL16} the authors of the present paper constructed a 2-arc-transitive non-normal Cayley graph on $\A_{39}$.
More results about normality of Cayley graphs we refer the reader to a survey paper in \cite{FLX08}.

Examples of connected arc-transitive non-normal cubic, tetravalent and pentavalent Cayley graphs
on nonabelian simple groups are very rare (the known examples are only the above mentioned
graphs on $\A_{47}$, $\M_{11}$, $\A_{35}$ and $\A_{39}$), we concentrate on the pentavalent case in this paper.
In particular, we construct a connected 3-arc-transitive non-normal pentavalent Cayley graph on $\A_{79}$ in Construction \ref{A79}.
It is shown in \cite{Li05} that 3-arc-transitive Cayley graphs of any given valency are rare.

The aim of this paper is to investigate the normality problem of the connected arc-transitive pentavalent Cayley graphs
with soluble vertex stabilizer on finite nonabelian simple groups.
Our main result is the following theorem.
\begin{thm}\label{thm1}
Let $G$ be a finite nonabelian simple group, and let $\Ga=\Cay(G,S)$ be an arc-transitive pentavalent Cayley graph on $G$.
Then the following statements hold.
\begin{itemize}
\item[(1)] Either $\Ga$ is a normal Cayley graph or $G=\A_{39}$ or $\A_{79}$. Further,
\item[(2)] there exist connected arc-transitive pentavalent non-normal Cayley graphs for $(G,\Aut\Ga)=(\A_{39},\A_{40})$ or $(\A_{79},\A_{80})$.
\end{itemize}
\end{thm}
\begin{rmk}\rm
(a) The connected 2-arc-transitive non-normal pentavalent Cayley graph on $\A_{39}$ in part (2) was constructed
by the authors of the present paper in \cite[Construction 3.1]{LL16}.

(b) The connected 3-arc-transitive non-normal pentavalent Cayley graph on $\A_{79}$ in part (2) constructed
in Section 4 is the first known example in 3-arc-transitive case.
\end{rmk}

\section{Preliminaries}
We give some necessary preliminary results in this section.
The first one is a property of the Fitting subgroup,
see \cite[P. 30, Corollary]{Suzuki}.

\begin{lemma}\label{Fitting-sg}
Let $F$ be the Fitting subgroup of a group $G$. If $G$ is soluble, then $F\ne 1$ and
the centralizer $\C_G(F)\le F$.
\end{lemma}

The next lemma is about primitive permutation groups of degree less than 80, refer to \cite{Ron05}.
\begin{lemma}\label{sg}
Let $T$ be a primitive permutation group on $\Omega$ and let $K$ be the stabilizer of a point $w\in\Omega$.
If $T$ is a nonabelian simple group and $|\Omega|$ divides 80, then $(T,K,|\Omega|)$ is one of the following Table \ref{table1}.
\begin{table}[ht]
 \[\begin{array}{|c|c|c|c|c|c|c|c|} \hline
T & \A_8 & \A_{10} & \A_{16} & \A_{20}&\PSL(4,3)&\A_{40}&\A_{80} \\ \hline
K &  \A_7 &  \A_9& \A_{15}& \A_{19}&\ZZ_3^3:\PSL(3,3)&\A_{39}&\A_{79}  \\ \hline
       |\Omega|&8&10&16&20&40&40&80\\ \hline
\end{array}\]
\caption{Primitive permutation groups of degree less than 80}\label{table1}
\end{table}
\end{lemma}
Simple groups which have subgroups of index dividing $2^5{\cdot}3^2$ are given in the following lemma, refer to \cite[Lemma 2.4]{FangLiXu-01}.
\begin{lemma}\label{simple-groups}
Let $T$ be a non-abelian simple group which has a subgroup $L$
of index dividing $2^5{\cdot}3^2$. Then $T$, $L$ and $n:=|T : L|$ are given in the following Table \ref{table-2}.
\begin{table}[ht]
 \[\begin{array}{ccccc}
      \hline
      T &L &n& {\rm Remark} \\
      \hline
       \A_n &\A_{n-1}&n  & n\ |\ 2^5{\cdot}3^2\\
        \ \M_{11} &\PSL(2,11)&  12 & \\
              \  \M_{12} &\M_{11}&  12 & \\
                \    \M_{24} &\M_{23}&  24 & \\  \hline
        \end{array}\]
\caption{Simple groups with having subgroups of index dividing $2^5{\cdot}3^2$}\label{table-2}
\end{table}
\end{lemma}

We next introduce the definition of coset graph. Let $G$ be a finite group and
let $H$ be a core-free subgroup of $G$.
Define the {\em coset graph} $\Cos(G,H,g)$ of $G$ with respect to $H$ as the graph with vertex set $[G:H]$
such that $Hx$, $Hy$ are adjacent if and only if $yx^{-1}\in HgH$. The following lemma about coset
graphs are well known and the proof of the lemma follows from the definition of coset graphs.
\begin{lemma}\label{lemma-coset}
Using notation as above. Let $\val\Ga$ be the valency of $\Ga$.
Then the coset graph $\Ga=\Cos(G,H,g)$ is $G$-arc-transitive graph and
\begin{itemize}
\item[(1)] $\val\Ga=|H:H\cap H^g|$;
\item[(2)] $\Ga$ is connected if and only if $\l H,g \r=G$.
\item[(3)] If $G$ has a subgroup $R$ acting regularly on the vertices of $\Cos(G,H,g)$, then $\Cos(G,H,g)\cong\Cay(R,S)$, where $S=R\cap HgH$.
\end{itemize}
Conversely, each $G$-arc-transitive graph $\Sigma$ is isomorphic to the coset graph
$\Cos(G,G_v,g)$, where $g\in \N_{G}(G_{vw})$ is a
$2$-element such that $g^2\in G_v$, and $v\in V\Sigma$, $w\in \Sigma(v)$.
\end{lemma}

For a graph $\Ga$ and a vertex-transitive subgroup $X\le\Aut\Ga$. Let $N$ be an intransitive normal subgroup of $X$ on
$V\Ga$. Denote $V_N$ the set of $N$-orbits in $V\Ga$. The {\em normal quotient graph} $\Ga_N$ defined as the
graph with vertex set $V_N$ and two $N$-orbits $B,C\in V_N$ are adjacent in $\Ga_N$ if and only if some
vertex of $B$ is adjacent in $\Ga$ to some vertex of $C$.
By \cite[Theorem 9]{Lorimer}, we have the following lemma.
\begin{lemma}\label{lor}
Let $\Ga$ be an arc-transitive graph of prime valency $p>2$ and let
$X$ be an arc-transitive subgroup of $\Aut\Ga$.
If a normal subgroup $N$ of $X$ has more than two orbits
on $V\Ga$, then $\Ga_N$ is an $X/N$-arc-transitive graph of valency
$p$ and $N$ is semiregular on $V\Ga$.
\end{lemma}

The following lemma is about the stabilizers of arc-transitive pentavalent graphs,
refer to \cite{Guo-1,Zhou-1}.

\begin{lemma}\label{stabilizer}
Let $\Ga$ be a pentavalent $(G,s)$-transitive graph, where $G\le\Aut\Ga$ and $s\ge 1$.
Let $\a\in V\Ga$. Then one of the following holds, where $\D_{10}$, $\D_{20}$ and $\F_{20}$
denote the dihedral groups of order $10$, $20$, and the Frobenius
group of order $20$, respectively.
\begin{itemize}
\item[(a)] If $G_{\a}$ is soluble, then $s\le 3$ and $|G_{\a}|\div 80$.
Further, the couple $(s,G_{\a})$ lies in the following Table \ref{table-3}.
\begin{table}[ht]
\[\begin{array}{|c|c|c|c|} \hline
s & 1 & 2 & 3 \\ \hline
G_{\a} &  \ZZ_5,~\D_{10},~\D_{20}& \F_{20},~\F_{20}\times\ZZ_2 & \F_{20}\times\ZZ_4  \\ \hline
\end{array}\]
\caption{The soluble case}\label{table-3}
\end{table}
\vskip0.2in

\item[(b)] If $G_{\a}$ is insoluble, then $2\le s\le 5$, and $|G_{\a}|\div 2^{9}\cdot 3^2\cdot 5$.
Further, the couple $(s,G_{\a})$ lies in the following Table \ref{table-4}.
\begin{table}[ht]
\[\begin{array}{|c|c|c|c|c|} \hline
s & 2 & 3 & 4 & 5 \\ \hline
G_{\a} &  \A_5,\S_5 &  \A_4\times\A_5,(\A_4\times\A_5){:}\ZZ_2,& \ASL(2,4),\AGL(2,4),& \ZZ_2^6{:}\GammaL(2,4)  \\
       & &  \S_4\times\S_5 & \ASigmaL(2,4),\AGammaL(2,4)&  \\ \hline
       |G_{\a}|&60,120&720,1440,2880&960,1920,2880,5760&23040\\ \hline
\end{array}\]
\caption{The insoluble case}\label{table-4}
\end{table}
\end{itemize}
\end{lemma}

\section{The proof of normal case }
%\noindent{\bf Proof of  Theorem \ref{thm1}}.
%=\ZZ_5$, $\D_{10}$, $\D_{20}$, $\F_{20}$, $\F_{20}\times\ZZ_2$ or $\F_{20}\times\ZZ_4$.
%If $\A_v=\ZZ_5$, $\D_{10}$ or $\D_{20}$, then by \cite[Theorem 5.4]{Zhou-1}, $G\unlhd\A$, Theorem \ref{thm1} holds.
%Hence we may assume that $\A_v=\F_{20}$, $\F_{20}\times\ZZ_2$ or $\F_{20}\times\ZZ_4$.
%If $G=\A_5$, then by \cite[Corollary 4.2]{ZF07},
%$G\unlhd\A$, Theorem \ref{thm1} holds. Hence we may also assume that $G\not=\A_5$.

Let $\Ga:=\Cay(G,S)$ be an arc-transitive pentavalent Cayley graph, where $G$ is a finite nonabelian simple group.
Let $\A:=\Aut\Ga$ and $\A_{v}$ be the stabilizer of $v$ in $\A$ where $v\in V\Ga$.
Assume that $\A_v$ is soluble. Then by Lemma \ref{stabilizer}, $|\A_v|$ divides 80.

The following lemma consider the case $\A$ has no nontrivial soluble normal subgroup.
\begin{lemma}\label{lemma1}
Assume that $\A$ has no nontrivial soluble normal subgroup. Then $G$ is either normal in $\A$ or $G=\A_{39}$ or $\A_{79}$.
\end{lemma}
\pf Let $N$ be a minimal normal subgroup of $\A$. Then $N=T^d$, where $d\ge1$ and $T$ is a nonabelian simple group.

Assume that $G$ is not normal in $\A$. Then since $N\cap G\unlhd G$ and $G$ is a nonabelian simple group,
$N\cap G=1$ or $G$. If $N\cap G=1$ then since $\A=G\A_v$, $|N|\div|\A_v|\div 80$,
which is a contradiction since $N$ is insoluble. Hence $N\cap G=G$, $G\le N$.
If $G=N$, then $G\unlhd\A$, a contradiction to the assumption. Thus $G<N$.
Assume that $d\ge2$. Then $N=T_1\times T_2\times\ldots\times T_d$ where $d\ge2$ and $T_i\cong T$ is a nonabelian simple group.
Since $T_1\cap G\unlhd G$, we have $T_1\cap G=1$ or $G$. If $T_1\cap G=1$, then $|T_1|\div|\A_v|\div80$, a contradiction.
If $T_1\cap G=G$, then $G\le T_1$. It follows that $|T_2|\div|\A_v|\div80$, which is also a contradiction.
Thus, $d=1$ and $N=T$ is a nonabelian simple group.
Then $T=GT_v$, $T_v \not= 1$ and $|T_v|$ divides 80.
Since $T$ has the proper subgroup $G$ with index dividing 80, we can take a maximal proper subgroup
$K$ of $T$ which contains $G$ as a subgroup. Let $\Omega=[T : K]$. Then $|\Omega|$ divides 80 and $T$ has
a primitive permutation representation on $\Omega$, of degree $n :=|\Omega|$. Since $T$ is simple, this
representation is faithful and thus $T$ is a primitive permutation group of degree $n$. Note
that $K$ is the stabilizer of a point $w\in \Omega$, that is, $K = T_w$. Since $T$ is nonabelian simple,
$n>4$. Consequently, by Lemma \ref{sg}, we have $(T,K,|\Omega|)$ is listed in Table \ref{table1}.

Assume that $(T,K,|\Omega|)=(\PSL(4,3),\ZZ_3^3:\PSL(3,3),40)$. Then since $G\le K$ and $G$ is a nonabelian simple group,
we have that $G$ is a proper subgroup of $K$. Since $|T:G|\div80$ and $|\Omega|=40$, we have $|K:G|=2$.
It is easy to see that a subgroup of $K$ with index 2 can not be a nonabelian simple group, which is a contradiction.

Assume that $(T,K,|\Omega|)=(\A_{8},\A_{7},8)$ or $(\A_{16},\A_{15},16)$.
Then since $G\le K$, $G$ is a nonabelian simple group and $|T:G|\le80$, we have $G=K$.
Since $\Ga$ is connected, $T\unlhd\A$ and $T_v\not=1$,
we have $1\not=T_v^{\Ga(v)}\unlhd\A_v^{\Ga(v)}$. Since $\Ga$ is $\A$-arc-transitive of valency 5, it follows that
$\A_v^{\Ga(v)}$ is primitive on $\Ga(v)$ and so $T_v^{\Ga(v)}$ is transitive on $\Ga(v)$, $5\div|T_v|=|T:G|=|\Omega|$, a contradiction.
Finally, by \cite[Theorem 5.4]{Zhou-1}, $(T,K,|\Omega|)\not=(\A_{20},\A_{19},20)$ or $(\A_{10},\A_{9},10)$.
Thus, we have $G=\A_{39}$ or $\A_{79}$, the lemma holds.
\qed
The following lemma consider the case $\A$ has a nontrivial soluble normal subgroup.
\begin{lemma}\label{lemma2}
Assume that $\A$ has a nontrivial soluble normal subgroup. Then $G$ is either normal in $\A$ or $G=\A_{39}$.
\end{lemma}
\pf Let $M$ be the largest soluble normal subgroup of $\A$. Then $M {\sf\,char\,}\A$.
Since $\A$ has a nontrivial soluble normal subgroup, we have $M\not=1$.
Since $M\cap G\unlhd G$ and $G$ is simple, we have $M\cap G = 1$ and so $|M|\div|\A_v|\div80$.
Since $|V\Ga|=|G|$ contains at least three prime factors, it follows that $M$ has more than two orbits on $V\Ga$.
By Lemma \ref{lor}, $M$ is semi-regular on $V\Ga$.

Let $\bar\A = \A/M$ and let $\bar\Ga=\Ga_M$. Then by Lemma \ref{lor}, $\bar\Ga$ is $\bar\A$-arc-transitive.
Let $\bar N$ be a minimal normal subgroup of $\bar\A$ and let $N$ be the full preimage of $\bar N$ under $\A\rightarrow\A/M$.
By the maximality of $M$, $\bar N$ is insoluble. Then $\bar N=T_1\times T_2\times \ldots\times T_d=T^d$, where $T$ is a nonabelian simple group and $d\ge1$.

We first show that $d=1$. Let $\bar G=GM/M$. Then $\bar G\cong G$ is a nonabelian simple group.
Since $\bar N\cap \bar G\unlhd\bar G$, we have $\bar N\cap\bar G=1$ or $\bar G$.
If $\bar N\cap\bar G=1$, then $|\bar N|$ divides 80, which is a contradiction since $\bar N$ is insoluble.
Hence $\bar G\le\bar N$. Since $\bar G$ is simple, $|\bar G|$ must divide the order of some composition factor
of $\bar N$, that is, $|\bar G|\div|T_1|$. If $d\ge2$ then $|T_2|$ divides $|\bar N : \bar G|$ which divides $|\bar\A_{\bar v}|$, which is not
possible since $\bar\A_{\bar v}$ is a $\{2,5\}$-group and $T_2$ is simple, where $\bar v\in V\bar\Ga$.
Therefore, $d=1$ and $\bar N$ is a nonabelian simple group. This argument also proves
that $\bar N$ is the unique insoluble minimal normal subgroup of $\bar \A$. Thus $\bar N{\sf\,char\,}\A$ and $N{\sf\,char\,}\A$.

Assume first that $\bar G =\bar N$. Then $N = M:G$. If $G$ centralizes $M$ then $N = M\times G$,
and therefore $G{\sf\,char\,} N {\sf\,char\,}\A$, a contradiction. Thus $G$ does not centralize $M$.
It follows that $\Aut(M)$ is insoluble.

Let $F$ be the Fitting subgroup of $M$.
By Lemma~\ref{Fitting-sg}, $F\ne 1$ and $\C_{M}(F)\le F$. Since $|M|$ divides 80,
we have $F=\O_2(M)\times\O_{5}(M)$,
where $\O_2(M)$, $\O_{5}(M)$ denote the largest normal $2$-, $5$-subgroups of $M$, respectively.
Clearly $|\O_2(M)|\div16$ and $|\O_5(M)|\div5$. Assume that $|\O_2(M)|\le2$.
Then $F$ is abelian and $F=\C_M(F)$. Since $M/\C_M(F)\lesssim\Aut(F)$, we have $M\lesssim F.\Aut(F)$.
If $|\O_2(M)|=1$, then since $\O_5(M)\le\ZZ_5$, we have $F\le\ZZ_5$.
It follows that $M\lesssim \ZZ_5.\ZZ_4$, and so $\Aut(M)$ is soluble, a contradiction.
If $|\O_2(M)|=2$, then $F\le\ZZ_{10}$.
Thus, $M\lesssim \ZZ_{10}.\ZZ_{4}$. A computation by Magma \cite{Magma}, $\Aut(M)$ is soluble, a contradiction.
Hence $|\O_2(M)|\ge4$.

Let $R =\O_2(M)$. Then $R{\sf\,char\,} M {\sf\,char\,}\A$. Let $B=RG$. We claim that $B$ is
not normal in $\A$. Suppose to the contrary that $B$ is normal in $\A$.
Then $B_v^{\Ga(v)}\unlhd\A_v^{\Ga(v)}$. Note that $\A_v^{\Ga(v)}$ is primitive on $\Ga(v)$.
Since $B>G$, $B_v\not= 1$ and therefore $B_v^{\Ga(v)}\not= 1$ is transitive on $\Ga(v)$.
Thus 5 divides $|B_v| = |B : G| = |R|$, a contradiction. So $B$ is not normal in $\A$ as claimed.

Assume that $G$ does not centralize $R$. Since $R\unlhd B$, we have that $B/\C_B(R)$ is isomorphic
to a subgroup of $\Aut(R)$. Since $G$ is nonabelian simple, it follows that $|R|\ge8$ and so
$|M/R|\le80/8=10$. Therefore $\Aut(M/R)$ is soluble.
Since $N/R=(M/R):(RG/R)$ and $RG/R\cong G$ is nonabelian simple,
we have $RG/R$ centralizes $M/R$, and so $N/R=(M/R)\times(RG/R)$.
It follows that $RG/R{\sf\,char\,} N/R$, $RG{\sf\,char\,} N\unlhd\A$,
which is a contradiction to the conclusion in the previous paragraph.

Thus $G$ centralizes $R$. Since $G$ does not centralize $M$, $R\not= M$.
Since $|R|\ge4$, it follows that $|M/R|\div20$, and so $\Aut(M/R)$ is soluble.
Similar arguments to the previous paragraph lead to $RG\unlhd\A$, a contradiction.

Thus $\bar G\not=\bar N$, $\bar G$ is a proper subgroup of $\bar N$ and $|\bar N:\bar G|$ divides 40.
Let $\bar K$ be a maximal proper subgroup of $\bar N$ which contains $\bar G$ as a subgroup and let $\bar\Omega=[\bar N : \bar K]$.
Then by Lemma \ref{sg}, we have $(\bar N,\bar K,|\bar\Omega|)$ is listed in Table \ref{table1}.
Since $|\bar\Omega|$ divides 40, we have $(\bar N,\bar K,|\bar\Omega|)\not=(\A_{80},\A_{79},80)$.
If $(\bar N,\bar K,|\bar\Omega|)=(\PSL(4,3),\ZZ_3^3:\PSL(3,3),40)$, then since $|\bar N:\bar G|\div40$ and $|\bar\Omega|=40$,
we have $\bar G=\bar K$, which is a contradiction since $\bar G$ is a nonabelian simple group.
Further, by the proof in \cite[Theorem 5.4]{Zhou-1}, $(\bar N,\bar K,|\bar\Omega|)\not=(\A_{20},\A_{19},20)$ or $(\A_{10},\A_{9},10)$.
Thus, to complete the proof of the lemma, we only need to exclude the cases where
$(\bar N,\bar K,|\bar\Omega|)=(\A_{16},\A_{15},16)$ or $(\A_{8},\A_{7},8)$.

Suppose that $(\bar N,\bar K,|\bar\Omega|)=(\A_{16},\A_{15},16)$ or $(\A_{8},\A_{7},8)$.
We claim that there exists $L\unlhd N$ such that 5 divides $|N : L|$.
Let $\bar v$ be a vertex in $V\bar\Ga$. Then $|\bar N_{\bar v}:\bar G_{\bar v}| = |\bar N : \bar G| = 16$ or 8.
Since $\bar\Ga$ is a pentavalent $\bar \A$-arc-transitive graph, we have $\bar\A_{\bar v}^{\bar\Ga(\bar v)}$ is
primitive on $\bar\Ga(\bar v)$. Since $1\not=\bar N_{\bar v}^{\bar\Ga(\bar v)}\unlhd\bar\A_{\bar v}^{\bar\Ga(\bar v)}$,
$\bar N_{\bar v}^{\bar\Ga(\bar v)}$ is transitive on $\bar\Ga(\bar v)$ and so $5$ divides $|\bar N_{\bar v}|$.
Therefore, $5$ divides $|\bar G_{\bar v}|$. Since $G$ is regular on $V\Ga$, it follows that $|M|=|\bar G_{\bar v}|$.
If $(\bar N,\bar G,|\bar\Omega|)=(\A_{16},\A_{15},16)$, then $|M|=5$.
Consequently, $N\cong\ZZ_5.\A_{16}$, an extension of $\ZZ_5$ by $\A_{16}$.
By Atlas \cite{ATLAS}, the Schur multiplier of $\A_{16}$ equals
$\ZZ_2$, and therefore, $N\cong\ZZ_5\times\A_{16}$. So there exists $L\cong\A_{16}$ as claimed.
If $(\bar N,\bar G,|\bar\Omega|)=(\A_{8},\A_{7},8)$, then $|M|=5$ or 10.
If $|M|=5$, then arguing as for the case $\bar N=\A_{16}$, there exists $L\cong\A_8$ as claimed.
Thus we suppose that $|M|=10$. Then $M\cong\ZZ_{10}$ or $\D_{10}$, and we can conclude that
$N\cong M_5\times(\ZZ_2.\A_8)$ or $(M_5\times\A_8).\ZZ_2$, where $M_5\cong\ZZ_5$ is a Sylow 5-subgroup of $M$.
Thus $L\cong\ZZ_2.\A_8$ or $\A_8$ exists as claimed. Since $|L|\not=|G|$,
$L$ is not regular on $V\Ga$. Now $1\not=L_v\unlhd N_v$, and since $N_v$ is transitive (and so primitive) on $\Ga(v)$,
$L_v$ is transitive on $\Ga(v)$ and so 5 divides $|L_v|$. Consequently, $5^2$ divides $|N:L|\cdot|L:G|=|N:G|=|N_v|$,
which is a contradiction to Lemma \ref{stabilizer}. Therefore, $G\cong\bar G\not\cong\A_{15}$ or $\A_7$,
and so $G\cong\A_{39}$. This completes the proof of the lemma.
\qed

\section{A 3-arc-transitive non-normal pentavalent Cayley graph on $\A_{79}$}
By \cite[Theorem 1.1]{LL16}, there exists a connected 2-arc-transitive non-normal pentavalent Cayley graph on $\A_{39}$ with full automorphism group $\A_{40}$.
In this section, we will construct a connected 3-arc-transitive non-normal pentavalent Cayley graph on $\A_{79}$
and prove its full automorphism group isomorphic to $\A_{80}$.

\begin{construct}\label{A79}
Let $G:={\sf Alt}(\{2,3,\ldots,80\})=\A_{79}$ and let $H=\l a,b,c\r<X:={\sf Alt}(\{1,2,\ldots,80\})=\A_{80}$, where
\[
\begin{array}{lllll}
a{=}&(1\ 16\ 11\ 6)(2\ 17\ 12\ 7)(3\ 18\ 13\ 8)(4\ 19\ 14\ 9)(5\ 20\ 15\ 10) \\
&(21\ 36\ 31\ 26)(22\ 37\ 32\ 27)(23\ 38\ 33\ 28)(24\ 39\ 34\ 29)(25\ 40\ 35\ 30) \\
&(41\ 56\ 51\ 46)(42\ 57\ 52\ 47)(43\ 58\ 53\ 48)(44\ 59\ 54\ 49)(45\ 60\ 55\ 50) \\
&(61\ 76\ 71\ 66)(62\ 77\ 72\ 67)(63\ 78\ 73\ 68)(64\ 79\ 74\ 69)(65\ 80\ 75\ 70), \\
b{=}&(1\ 46\ 77\ 35)(2\ 43\ 66\ 28)(3\ 60\ 75\ 21)(4\ 57\ 64\ 34)(5\ 54\ 73\ 27) \\
&(6\ 51\ 62\ 40)(7\ 48\ 71\ 33)(8\ 45\ 80\ 26)(9\ 42\ 69\ 39)(10\ 59\ 78\ 32) \\
&(11\ 56\ 67\ 25)(12\ 53\ 76\ 38)(13\ 50\ 65\ 31)(14\ 47\ 74\ 24)(15\ 44\ 63\ 37) \\
&(16\ 41\ 72\ 30)(17\ 58\ 61\ 23)(18\ 55\ 70\ 36)(19\ 52\ 79\ 29)(20\ 49\ 68\ 22),\\
c{=}&(1\ 17\ 13\ 9\ 5)(2\ 18\ 14\ 10\ 6)(3\ 19\ 15\ 11\ 7)(4\ 20\ 16\ 12\ 8) \\
&(21\ 37\ 33\ 29\ 25)(22\ 38\ 34\ 30\ 26)(23\ 39\ 35\ 31\ 27)(24\ 40\ 36\ 32\ 28) \\
&(41\ 57\ 53\ 49\ 45)(42\ 58\ 54\ 50\ 46)(43\ 59\ 55\ 51\ 47)(44\ 60\ 56\ 52\ 48) \\
&(61\ 77\ 73\ 69\ 65)(62\ 78\ 74\ 70\ 66)(63\ 79\ 75\ 71\ 67)(64\ 80\ 76\ 72\ 68).
\end{array}
\]
Take $x_1\in G$ as follows:
\[
\begin{array}{ll}
x_1{=}&(2\ 22)(3\ 29)(4\ 36)(5\ 23)(6\ 35)(7\ 68)(8\ 79)(9\ 70)(10\ 61)(11\ 77) \\
&(12\ 49)(13\ 52)(14\ 55)(15\ 58)(16\ 46)(17\ 20)(18\ 19)(21\ 34)(24\ 60) \\
&(25\ 62)(26\ 64)(27\ 28)(30\ 51)(31\ 57)(32\ 66)(33\ 73)(37\ 43)(38\ 54) \\
&(39\ 75)(42\ 65)(44\ 53)(45\ 74)(47\ 50)(48\ 63)(56\ 72)(59\ 76)(69\ 80)(71\ 78).
\end{array}
\]
Define $\Sigma=\Cos(X,H,x_1)$.
\end{construct}
\begin{lemma}\label{lem-A80}
The graph $\Sigma=\Cos(X,H,x_1)$ in Construction \ref{A79} is connected,
3-arc-transitive and isomorphic to the non-normal Cayley graph $\Cay(G,S)$ of $G$, determined
by $S=\{x_1,x_2,x_2^{-1},x_3,x_3^{-1}\}$ with
\[\begin{array}{llll}
x_2{=}&(2\ 5\ 40\ 79\ 75\ 49\ 67\ 36\ 10\ 74\ 37\ 8\ 72\ 62\ 14\ 56\ 18\ 4\ 33\ 70\ 64\ 52\ 34\ \\
&77\ 43\ 69\ 11\ 65\ 30\ 7\ 39\ 58\ 35\ 25\ 44\ 21\ 24\ 9\ 63\ 29\ 31\ 55\ 22\ 38\ 47\ 20) \\
&(3\ 26\ 66\ 46\ 78\ 15\  59\ 19)(6\ 27\ 68\ 71\ 13\ 53\ 17\ 50\ 16\ 42\ 45\ 57\ 76\ 23) \\
&(12\ 61\ 73\ 60\ 80)(28\ 32\ 41)(48\ 54), \\
x_3{=}&(2\ 20)(3\ 4\ 25\ 8\ 38\ 7\ 26\ 43\ 44\ 68\ 35\ 76\ 48\ 19)(5\ 32\ 54\ 16\ 58\ 75\ 36\ 24) \\
&(6\ 39\ 46\ 56\ 21\ 78\ 42\ 27\ 65\ 51\ 15\ 55\ 47\ 77\ 22\ 23\ 37\ 9\ 71\ 63\ 13\ 80\ 74\ \\
&59\ 34\ 69\ 70\ 10\ 62\ 72\ 12\ 64\ 28\ 67\ 45\ 40\ 57\ 17\ 41\ 79\ 14\ 52\ 50\ 66\ 29\ 31) \\
&(11\ 73\ 61)(18\ 49\ 53\ 33\ 60).
\end{array}
\]
\end{lemma}
\pf Let $\Delta=\{1,2,\ldots,80\}$ and consider the natural action of $X$ on $\Delta$.
By Magma \cite{Magma}, $\l H,x_1 \r=X$, and so $\Sigma$ is connected by Lemma \ref{lemma-coset} $(2)$.
Since $c^b=c^2$ and $a$ centralizes $\l b,c\r$, it follows that
$H=\l a,b,c\r=\l b,c \r{\times}\l a \r\cong (\ZZ_5{:}\ZZ_4){\times}\ZZ_4$. Furthermore, it is
easy to see that $H$ is transitive on $\Delta$ and so is regular on $\Delta$. Hence
$X$ has a factorization $X=GH=HG$ with $G\cap H=1$. Therefore, $\Sigma$ is
isomorphic to a Cayley graph of $G=\A_{79}$.
Further computation shows that $\frac{|H|}{|H\cap H^{x_1}|}=5$.
By Lemma \ref{lemma-coset} $(1)$ we have that $\Sigma$ is pentavalent.
Since $H\cong (\ZZ_5{:}\ZZ_4){\times}\ZZ_4$, we have $\Sigma$ is 3-arc-transitive by Lemma \ref{stabilizer}.
Since $X$ is simple, $G$ is not normal in $X\le\Aut\Sigma$. Thus $\Sigma$ is non-normal.
Let $x_2$, $x_3$ and $S$ define as in this lemma.
Computation shows that $G\cap (Hx_1H)=S$.
Then $\Sigma\cong \Cay(G,S)$ by Lemma \ref{lemma-coset} $(3)$.
This completes the proof of the lemma.
\qed
In the next lemma we show that $\Aut\Sigma\cong\A_{80}$.
This will therefore complete the proof of Theorem \ref{thm1}.
\begin{lemma}\label{lemma4}
The full automorphism group $\Aut\Sigma$ of the graph $\Sigma=\Cos(X,H,x_1)$ in Construction \ref{A79} is isomorphic to $\A_{80}$.
\end{lemma}
\pf Let $\A=\Aut\Sigma$. Assume first that $\A$ is quasiprimitive on $V\Sigma$.
Let $N$ be a minimal normal subgroup of $\A$. Then $N$ acts transitively on
$V\Sigma$, and so $N$ is insoluble.
Then $N = T^d$ with $T$ a nonabelian simple group and $d\ge1$.
Let $p$ be the largest prime factor of $|\A_{79}|$. Then $p>5$ and $p^2\not\big|\,|\A_{79}|$. Since $N$ is transitive on
$V\Sigma$ and $|V\Sigma|=|\A_{79}|$, we have $p\div|N|$.
Suppose that $d\ge2$. Then $p^d\div|N|$. However, by Lemma \ref{stabilizer},
$|\A_v|\div 2^9\cdot3^2\cdot5$, and so $p^d\div|N|\div|\A|\div2^{9}\cdot3^2\cdot5\cdot|\A_{79}|$,
a contradiction. Hence $d = 1$ and $N= T\unlhd \A$. Let $C=\C_{\A}(T)$. Then $C\unlhd \A$
and $CT = C{\times}T$. If $C\not=1$, then $C$ is transitive on $V\Sigma$ as
$\A$ is quasiprimitive on $V\Sigma$. It follows that $p\div|C|$.
Therefore, $p^2\div|CT|\div |\A|$, a contradiction.
Hence $C = 1$ and $\A\le\Aut(T)$ is almost simple.

Since $T\cap X \unlhd X\cong\A_{80}$, it follows that $T\cap X=1$ or $\A_{80}$.
If $T\cap X=1$, then $|T|\div\frac{|\A|}{|X|}\div2^5\cdot3^2$, and so $T$ is soluble, a contradiction.
Thus $T\cap X=X$, and so $X\le T$.
It follows that $|T:X|\div|\A:X|\div2^5\cdot3^2$.
By Lemma \ref{simple-groups} we can conclude that $T=X\cong\A_{80}$.
Thus $\A\le \Aut(T)\cong\S_{80}$. If $\A\cong\S_{80}$, then $|\A_v|=\frac{|\A|}{|G|}=160$,
a contradiction with Lemma \ref{stabilizer}.
Hence $\A\cong\A_{80}$.

Now assume that $\A$ is not quasiprimitive on $V\Sigma$.
Let $M$ be a minimal normal subgroup of $\A$ which is not transitive on $V\Sigma$.
Then $M\cap X \unlhd X$. It follows that $M\cap X =1$ or $\A_{80}$. If
$M\cap X =\A_{80}$, then $X\le M$, and so $M$ is transitive on $V\Sigma$, a contradiction.
If $M\cap X =1$, then $|M|$ divides $\frac{|\A|}{|X|}\div2^5\cdot3^2$.
Thus $M\cong\ZZ_2^r$ or $\ZZ_3^l$, where $1\le r\le 5$ and $1\le
l\le2$. Let $L=MX$.
Then $L=M{:}X$ and $L/\C_L(M)\lesssim \Aut(M)\cong \GL(r,2)$ or $\GL(l,3)$. Note
that $M\le \C_L(M)$. If $M=\C_L(M)$, then $L/\C_L(M)=L/M\cong
X \cong\A_{80}\lesssim \GL(r,2)$ or $\GL(l,3)$. However,
$\GL(2,r)$ or $\GL(3,l)$ has no subgroup
isomorphic to $\A_{80}$ for $1\le r\le 5$ and $1\le l\le 2$. Hence we
have $M<\C_L(M)$ and $1\not=\C_L(M)/M\unlhd L/M\cong\A_{80}$. It
follows that $\A_{80}=\C_L(M)/M$, that is,
$X$ centralizes $M$. Hence $L=M\times X$. Then
$L_v/X_v\cong L/X\cong M$. It implies
that $L_v$ is soluble and $L_v\cong X_v.M$. Since $|X_v|=80$,
we have $|L_v|>80$, a contradiction with Lemma \ref{stabilizer}.
This completes the proof of the lemma. \qed

\end{document}